\documentclass{amsart}
\pdfoutput=1
\usepackage{graphicx}
\newtheorem{thm}{Theorem}[section]
\newtheorem{cor}[thm]{Corollary}
\newtheorem{lem}[thm]{Lemma}
\newtheorem{prop}[thm]{Proposition}
\theoremstyle{definition}
\newtheorem{defn}[thm]{Definition}
\theoremstyle{remark}

\numberwithin{equation}{section}

\begin{document}
	
	\title[]{Vilenkin-Fourier series in variable Lebesgue spaces}%
	
	\author{Daviti Adamadze and Tengiz Kopaliani} 
	
\address{
	Faculty of Exact and Natural Sciences\\
	Javakhishvili Tbilisi State
	University\\
	13, University St., Tbilisi, 0143, Georgia
}
\email{daviti.adamadze2013@ens.tsu.edu.ge}

\address{
			Faculty of Exact and Natural Sciences\\
			Javakhishvili Tbilisi State
			University\\
			13, University St., Tbilisi, 0143, Georgia}
		\email{tengiz.kopaliani@tsu.ge}
	\date{29.01.2022}
\thanks{This work was supported by Shota Rustaveli National Science Foundation of Georgia FR-21-12353}


\begin{abstract}
	Let $S_{n}f$ denote the $n$th partial sum of the Vilenkin-Fourier series of a function $f \in L^{1}(G)$. For $1 < p_{-} \leq p_{+} < \infty$, we characterize all exponents $p(\cdot)$ for which the convergence of $S_{n}f$ to $f$ in $L^{p(\cdot)}(G)$ holds whenever $f \in L^{p(\cdot)}(G)$.
\end{abstract}

	\maketitle
	\section{Introduction}
	
Let $\{p_{i}\}_{i\geq0}$ be a sequence of integers with $p_{i} \geq 2$. Define $G = \Pi_{i=0}^{\infty} \mathbb{Z}_{p_{i}}$ as the direct product
of cyclic groups of order $p_{i}$, and let $\mu$ be the Haar measure on $G$, normalized such that $\mu(G) = 1$. Each element of $G$ can be represented 
as a sequence $\{x_{i}\}$ with $0 \leq x_{i} < p_{i}$. Define $m_{0} = 1$ and $m_{k} = \Pi_{i=0}^{k-1} p_{i}$ for $k = 1,2,\dots$. 

There exists a well-known and natural measure-preserving identification between the group $G$ and the closed interval $[0,1]$. 
This identification is established by associating each sequence $\{x_{i}\} \in G$, where $0 \leq x_{i} < p_{i}$, with the point
$
\sum_{i=0}^{\infty} x_{i} m_{i+1}^{-1}.
$
Disregarding the countable set of $p_{i}$-rationals, this mapping is one-to-one, onto, and measure-preserving.

	For each $x=\{x_{i}\}\in G,$ define $\phi_{k}(x)=\exp(2\pi i x_{k}/p_{k}),$ $k=0,1,....$  The set $\{\psi_{n}\}$ of characters of $G$ consists of all finite product of $\phi_{k},$ which we enumerate in the following manner.
	Express each nonnegative integer $n$
	as a finite  sum $n=\sum_{i=0}^{\infty}\alpha_{k}m_{k}$, with $0\leq \alpha_{k}<p_{k},$ and define $\psi_{n}=\Pi_{i=0}^{\infty}\phi_{k}^{\alpha_{k}}.$
	The functions $\psi_{n}$ form a complete orthonormal system on $G.$  For the case $p_{i}=2,\,i=0,1,...,$  $G$ is the dyadic group, $\phi_{k}$ are Rademacher functions and $\psi_{n}$ are Walsh functions. In general, the system $\{\psi_{n}\}$ is a realization of the multiplicative Vilenkin system. In this paper, there is no restriction on the orders $\{p_{i}\}.$
	
	For $f\in L^{1}(G)$, let $S_{n}f,\,\,n=0,1,...$, be the $n$th partial sum of the Vilenkin-Fourier series of $f.$  When the orders $p_{i}$ of cyclic groups  are bounded  Watari \cite{Wat} showed that for $f\in L^{p}(G),\,\,1<p<\infty$,
	$$
	\lim_{n\rightarrow\infty}\int_{G}|S_{n}f-f|^{p}d\mu=0.
	$$
	Young \cite{You1}, Schipp \cite{Sc} and Simon \cite{Sim} showed independently that results concerning mean convergence of partial sums of the Vilenkin-Fourier series are still valid even if the orders $p_{i}$ are unbounded.
	
	Let $\{G_{k}\}$  be the sequence of subgroups of $G$ defined by
	\[
	G_{0}=G,\,\,G_{k}=\Pi_{i=0}^{k-1}\{0\}\times\Pi_{i=k}^{\infty}\mathbb{Z}_{p_{i}},\,\,\,k=1,2,....
	\]
	
	On the closed interval $[0,1]$, cosets of $G_{k}$  are intervals of the form  $[jm_{k}^{-1},(j+1)m_{k}^{-1}],$ $j=0,1,...,m_{k}-1.$
	By $\mathcal{F}$ we denote the set of generalized intervals. This set is
	the collection of all translations of  intervals $[0,jm_{k+1}^{-1}],\,k=0,1,...$ $j=1,...,p_{k}.$
	Note that a set $I$ belongs to $\mathcal{F}$  if (1) for some $x\in G$ and $k,$ $I\subset x+G_{k},$
	(ii) $I$ is a union of cosets of $G_{k+1,}$ and (iii) if we consider $ x+G_{k}$  as a circle, $I$ is an interval. Let $\mathcal{F}_{-1}=\{G\}$. For $k=0,1,...,$ let $\mathcal{F}_{k}$ be the collection of all $I\in \mathcal{F}$ such that $I$ is a proper subset of a coset of $G_{k}$, and is a union of cosets of $G_{k+1}.$ The collections $\mathcal{F}_{k}$ are disjoint, and $\mathcal{F}=\cup_{k=-1}^{\infty}\mathcal{F}_{k}.$ For $I\in\mathcal{F}$, we define the set $3I\in\mathcal{F}$ as follows. If $I=G$, let $3I=G$. For $I\in\mathcal{F}_{k}, k=0,1,...,$ there is $x\in G$ such that $I\subset x+G_k$. If $\mu(I)\geq\frac{\mu(G_k)}{3}$, let $3I=x+G_k$. If $\mu(I)<\frac{\mu(G_k)}{3}$, consider $x+G_k$ as a circle. Then $I$ is an interval in this circle. Define $3I\in\mathcal{F}_k$ to be the interval in this circle which contains $I$ at its center and has measure $\mu(3I)=3\mu(I)$. In all cases, for $I\in\mathcal{F}$, $\mu(3I)\leq 3\mu(I)$.

	We say that $w$ is a weight function on $G$ if $w$ is measurable and $0<w(x)<\infty$ a.e.  Gosselin \cite{Gos} (case $\sup_{i}p_{i}<\infty$ ) and  Young \cite{You2} (no restriction on the orders $p_{i}$) characterized all weight functions $w$ such that if $f\in L^{p}_{w}(G),\,\,1<p<\infty,$  $S_{n}f$ converges to $f$ in $ L^{p}_{w}(G)$. Here $L_{w}^{p}(G)$ denotes the space of measurable functions on $G$ such that $\|f\|_{p,w}=(\int_{G}|f|^{p}wd\mu)^{1/p}<\infty.$

	\begin{defn}(see \cite{You2})
		(i) We say that $w$ satisfies $A_{p}(G)$ condition, $1<p<\infty$, if
		\begin{equation}\label{Muk}
			[w]_{A_{p}}=\sup_{I\in\mathcal{F}} \left(\frac{1}{\mu(I)}\int_{I}w\,d\mu\right) \left(\frac{1}{\mu(I)}\int_{I}w^{-1/(p-1)}\,d\mu\right)^{p-1}<\infty.
		\end{equation}
		
		(ii) We say that $w$ satisfies $A_{1}(G)$ condition if 
		$$
		[w]_{A_{1}}=\sup_{I\in\mathcal{F}}  \frac{1}{\mu(I)}\int_{I}wd\mu\,(\mbox{essinf}_{I} w(x))^{-1}<\infty.
		$$
	\end{defn}
	
	For the case where the orders of cyclic groups are bounded, Gosselin \cite{Gos} defined $A_{p}(G)$
	condition, as the one where \eqref{Muk} condition holds
	for all $I$ that are cosets of $G_{k}, \,k=0,1,2,....$ For this case $A_{p}$ conditions, defined by  Young and  Gosselin,  are equivalent (see \cite{You2}).
	
	\begin{thm}(\cite{You2})
		Let $w$ be a weight function on $G.$ For $1<p<\infty$, the following  statements are equivalent:
		
		(i) $w\in A_{p}(G),$
		
		(ii) There is a constant $C$, depending only on $w$ and $p$, such that for every $f\in L^{p}_{w}(G)$, we have
		\[
		\int_{G}|S_{n}f|^{p}wd\mu\leq C\int_{G}|f|^{p}wd\mu,
		\]
		
		(iii) For every $f\in L^{p}_{w}(G)$, we have
		$$
		\lim_{n\rightarrow\infty}\int_{G}|S_{n}f-f|^{p}wd\mu=0.
		$$
	\end{thm}
	
	In this paper we characterize all exponents $p(\cdot)$ such that if $f\in L^{p(\cdot)}(G)$, then  partial sums $S_{n}f$ of the Vilenkin-Fourier series of  $f\in L^{p(\cdot)}(G)$ converge to $f$ with $L^{p(\cdot)}$-norm. Now we give a definition of variable Lebesgue space. Let $p(\cdot):G\rightarrow[1,\infty)$ be a measurable function.  The variable Lebesgue space $L^{p(\cdot)}(G)$ is the set of all measurable functions $f$ such that for some $\lambda>0,$
	$$
	\rho_{p(\cdot)}(f/\lambda)=\int_{G}(|f(x)|/\lambda)^{p(x)}d\mu<\infty.
	$$
	$L^{p(\cdot)}(G)$ is a Banach function space equipped with the Luxemburg norm
	$$
	\|f\|_{p(\cdot)}=\inf\{\lambda>0:\,\,\rho_{p(\cdot)}(f/\lambda)\leq1\}.
	$$
	
	We use the notations $p_{-}(I)=\mbox{essinf}_{x\in I}p(x)$ and $p_{+}(I)=\mbox{esssup}_{x \in I}p(x)$ where $I\subset G$. If $I=G$ we simply use the following notation $p_{-}, p_{+}$. The function $p'(\cdot)$  denotes the conjugate exponent function of $p(\cdot)$, i.e., $1/p(x)+1/p'(x)=1\,\,(x\in G).$ In this paper the constants $C, c$ are absolute constants and may be different in different contexts and $\chi_{A}$ denotes the characteristic function of set $A.$

	Very recently  the convergence of partial sums of the Walsh-Fourier series   in $L^{p(\cdot)}([0,1))$ space was investigated by Jiao et al. \cite{JWWZ}.  We denote by $C_{d}^{\log}$ the set of all functions $p(\cdot): [0,1)\rightarrow [1\infty)$, for which  there exists a positive constant $C$ such that
	$$
	|I|^{p_{-}(I)-p_{+}(I)}\leq C
	$$
	for all dyadic intervals $I=[k2^{-n}, (k+1)2^{-n})$ $(k,n\in \mathbb{N},\,0\leq k<2^{n})$, here $|I|$ denotes the Lebesgue measure of $I$.
	Note that this  condition may be interpreted as a dyadic version of log-H\"{o}lder continuity condition of $p(\cdot)$\,(or on dyadic group). The log-H\"{o}lder condition is a very common condition for solving various problems of harmonic analysis in $L^{p(\cdot)}(\mathbb{R}^{n})$ (see \cite{CUF}, \,\cite{DHHR}).
	
	\begin{thm}(\cite{JWWZ})
		Let $p(\cdot)\in C_{d}^{\log}$ with $1<p_{-}\leq p_{+}<\infty.$ If $f\in L^{p(\cdot)}([0,1))$, then for partial sums $S_n f$ of the Walsh-Fourier series of $f\in L^{p(\cdot)}([0,1))$ we have
		$$
		\sup_{n\in \mathbb{N}}\|S_{n}f\|_{p(\cdot)}\leq C \|f\|_{p(\cdot)}.
		$$
	\end{thm}
	
	Since  Walsh polynomials are dense in  $L^{p(\cdot)}([0,1))$, Theorem 1.3 implies that $S_{n}f$ converges to the original function in  $L^{p(\cdot)}([0,1))$-norm (for more details see \cite{JWWZ} and the recent book \cite{PTW}, chapter 9).
	
	In order to extend techniques and results of constant exponent case to the setting of variable Lebesgue spaces, a central problem is to
	determine conditions on an exponent $p(\cdot)$ under which the Hardy-Littlewood maximal operator is bounded on $L^{p(\cdot)}$ (see monographs Cruz-Uribe and Fiorenza \cite{CUF} and Diening et.al. \cite{DHHR}). We now define the Hardy-Littlewood maximal function that is appropriate for the study of Vilenkin-Fourier series. For $f\in L^{1}(G)$, let
	$$
	Mf(x)=\sup_{x\in I, I\in \mathcal{F}}\frac{1}{\mu(I)}\int_{I}|f|d\mu.
	$$
	
	This maximal function was introduced first by P. Simon in \cite{Sim1}.  He showed that the maximal operator is bounded in $L^{p}(G),\,1<p<\infty$ and is of weak type $(1,1).$
	Young \cite{You2} obtained the following analogue of Muckenhoupt's  theorem \cite{Muck}.
	\begin{thm}
		Let $w$ be a weight function on $G.$ For $1<p<\infty$, the following two statements are equivalent:
		
		(i) $w\in A_{p}(G),$
		
		(ii) There is a constant $C$, depending only on $w$ and $p$, such that for every $f\in L^{p}_{w}(G)$, we have
		\[
		\int_{G}(Mf)^{p}wd\mu\leq C\int_{G}|f|^{p}wd\mu.
		\]
		
		In case $p=1$ the following two statements are also equivalent:
		
		(iii) $w\in A_{1}(G),$
		
		(iv) There is constant $C,$ depending only on $w$, such that for every  $f\in L^{1}(G)$
		\[
		\int_{\{Mf>y\}}wd\mu\leq Cy^{-1}\int_{G}|f|wd\mu,\,\,\,y>0.
		\]
	\end{thm}

	\begin{defn}
		We say that the exponent $p(\cdot),\,1<p_{-}\leq p_{+}<\infty$  satisfies the condition $\mathcal{A}(G)$,  if there is a constant $C$ such that for every $I\in \mathcal{F}$,
		\begin{equation}\label{vMuk}
			\frac{1}{\mu(I)}\|\chi_{I}\|_{p(\cdot)}\|\chi_{I}\|_{p'(\cdot)}\leq C.
		\end{equation}
	\end{defn}
	
	The condition \eqref{vMuk} plays exactly the same role for averaging operators in variable Lebesgue spaces as the Muckenhoupt $A_{p}$ conditions for weighted Lebesgue spaces (see \cite{K1}, \cite{L1}, for Euclidian setting). We show that the $\mathcal{A}(G)$ condition is necessary and sufficient for the $L^{p(\cdot)}(G)$ boundedness of Hardy-Littlewood maximal function. One of the main result of the present paper is the following theorem.
	
	\begin{thm}\label{Muck}
		Assume for the exponent $p(\cdot)$ we have $1<p_{-}\leq p_{+}<\infty.$ Then the following two statements are equivalent:
		
		(i) $p(\cdot)\in \mathcal{A}(G),$
		
		(ii) There is a constant $C$, depending only on $p(\cdot)$ such that for every $f\in L^{p(\cdot)}(G)$, we have
		\[
		\|Mf\|_{p(\cdot)}\leq C\|f\|_{p(\cdot)}.
		\]
	\end{thm}
	
	By the symmetry of the definition, $p(\cdot)\in\mathcal{A}(G)$ if and only if  $p'(\cdot)\in\mathcal{A}(G)$ and from Theorem 1.6  we have  that, even though, $M$ is not a linear operator, the boundedness of $M$ implies the "dual" inequality.
	
	\begin{cor}
		Let for exponent $p(\cdot)$ we have $1<p_{-}\leq p_{+}<\infty.$ Then the maximal operator $M$ is bounded on $L^{p(\cdot)}(G)$ if and only if $M$ is bounded on $L^{p'(\cdot)}(G).$
	\end{cor}
	
	We prove the following theorem (in the Euclidean setting see \cite{CUF}, Theorem 4.37 and \cite{DHHR}, Theorem 5.7.2).
	
	\begin{thm} Let for the exponent $p(\cdot)$ we have $1<p_{-}\leq p_{+}<\infty.$
		Then the following  statements are equivalent:
		
		(i) Maximal operator $M$ is bounded on $L^{p(\cdot)}(G)$,
		
		(ii) There exists $r_{0}$, $0<r_{0}<1,$ such that if $r_{0}<r<1,$ then maximal operator $M$ is bounded on $L^{rp(\cdot)}(G)$.
	\end{thm}
	
	Hereafter, we will denote by $\mathcal{S}$  a family of pairs of non-negative,
	measurable functions.  Given $p$, $1\leq p<\infty$ if for some $w\in A_{p}(G)$ we write
	$$
	\int_{G}f(x)^{p}w(x)d\mu\leq C\int_{G}g(x)^{p}w(x)d\mu,\,\,\,\,\,(f,g)\in\mathcal{S},
	$$
	then we mean that this inequality holds for all pairs $(f,g)\in\mathcal{S}$ such that the left hand
	side is finite, and that the constant $C$ may depend on $p$ and $[w]_{A_{p}}$.
	If we write
	$$
	\|f\|_{p(\cdot)}\leq C_{p(\cdot)}\|g\|_{p(\cdot)},\,\,\,\,\,(f,g)\in\mathcal{S},
	$$
	then we mean that this inequality holds for all pairs $(f,g)\in\mathcal{S}$ such that the left-hand side is
	finite and the constant may depend on $p(\cdot)$.
	
	Using this convention we can state the Rubio
	de Francia extrapolation theorem in the following manner.
	
	\begin{thm}
		Suppose for some $p_{0}\geq1$ the family $\mathcal{S}$ is such that for all $w\in A_{1}(G)$
		$$
		\int_{G}f(x)^{p_{0}}w(x)d\mu\leq C\int_{G}g(x)^{p_{0}}w(x)d\mu,\,\,\,\,\,(f,g)\in\mathcal{S}.
		$$
		If for the exponent $p(\cdot)$, we have $p_{0}< p_{-}\leq p_{+}<\infty$ and the maximal operator M is bounded on $L^{(p(\cdot)/p_{0})'}(G)$, then
		$$
		\|f\|_{p(\cdot)}\leq C_{p(\cdot)}\|g\|_{p(\cdot)},\,\,\,\,\,(f,g)\in\mathcal{S}.
		$$
	\end{thm}
	
	Firstly, Theorem 1.9  was proved in \cite{CUFMP} (Theorem 1.3) for variable exponent  Lebesgue spaces on $\mathbb{R}^{n}$ and maximal
	operator $M$ defined on cubes (balls) in  $\mathbb{R}^{n}$, with sides parallel to the coordinate axes. In \cite{CUMP} the Rubio
	de Francia extrapolation theorem is proved  for general Function spaces, using $A_{1}$ weights and maximal
	operator $M$ defined by any Muckenhoupt basis (see Definition 3.1 in \cite{CUMP}). By Theorem 1.4  the set of
	generalized intervals $\mathcal{F}$ is a Muckenhoupt basis. Considering the following equality $\left(  L^{p(\cdot)}(G)\right) ^{1/p_0}=L^{p(\cdot)/p_0}(G)$, Theorem 1.9 is direct consequence of Theorem 4.6 from \cite{CUMP}.

	Now, we can formulate the main result of the present paper.
	\begin{thm}
		Let for exponent $p(\cdot)$ we have $1<p_{-}\leq p_{+}<\infty.$ Then the following  statements are equivalent:
		
		(i) $p(\cdot)\in \mathcal{A}(G),$
		
		(ii) There is a constant $C$, depending only on $p(\cdot)$, such that for  partial sums $S_n f$ of the Vilenkin-Fourier series of  $f\in L^{p(\cdot)}(G)$  we have
		$$
		\sup_{n\in \mathbb{N}}\|S_{n}f\|_{p(\cdot)}\leq C \|f\|_{p(\cdot)}.
		$$
		
		(iii) Partial sums $S_{n}f$ of the Vilenkin-Fourier series of  $f\in L^{p(\cdot)}(G)$ converge to the original function in $L^{p(\cdot)}$ space.
	\end{thm}

	\section{Preliminaries}

	The fundamental properties of $ A_{p}(G)$ weights were investigated by Gosselin \cite{Gos} and later by Young \cite{You2} (in this paper there is no restriction on the orders $p_{i}$). We formulate some  properties of these weights (see \cite{You2}).
	
	Note  that if $w\in A_{p}(G)$, then $L^{p}_{w}(G)\subset L^{1}(G)$.  We also mention  that if $w\in A_{p}(G),\,1\leq p<\infty$, and $p<q<\infty$
	then $w\in A_{q}(G)$. A important property of $A_{p}(G)$ weights is the reverse H\"{o}lder inequality.
	\begin{prop} (\cite{You2})
		Let $w\in A_{p}(G),\,1<p<\infty$. Then there exist $s>1$ and a constant $C$ such that for any $I\in\mathcal{F},$
		\[
		\left(\frac{1}{\mu(I)}\int_{I}w^{s}d\mu\right)^{1/s}\leq\frac{C}{\mu(I)}\int_{I}wd\mu.
		\]
	\end{prop}
	
	The following proposition is a consequence of the reverse H\"{o}lder inequality.
	
	\begin{prop}(\cite{You2})
		(i) Suppose $w\in A_{p}(G),1<p<\infty$. Then there exists $1<s<p$ such that $w\in A_{s}(G)$.
		(ii) Suppose $w\in A_{p}(G),1<p<\infty$, then $w\in A_{\infty}(G).$
	\end{prop}

	\begin{defn}(\cite{You2})
		Let $I_{0}\in \mathcal{F}$. We say that a weight $w$ (i.e. a nonnegative integrable function) satisfies $A_{\infty}(I_{0})$ condition if
		for any $\varepsilon\in(0,1)$  there exists $\delta\in(0,1)$ such that for any generalized interval $I\subset I_{0}$ and for any measurable subset $E\subset I$,
		$\mu(E)\leq\varepsilon \mu(I)$ implies $w(E)\leq \delta w(I)$ (for any measurable set $A$, $w(A)=\int_{A}wd\mu$ and $w_A=\frac{1}{\mu(A)}\int_{A}wd\mu$).
	\end{defn}
	
	It is well known fact that the class $A_{\infty}$ in Euclidian case can be defined in many equivalent ways. The most classical definition is due
	to Muckenhoupt \cite{Mu}. It is said that a locally integrable function $w:\mathbb{R}^{n}\rightarrow[0,\infty)$ is in $A_{\infty}$ class if for
	each $\varepsilon\in (0,1)$ there exists $\delta\in (0,1)$ such that $|E|\leq\varepsilon|Q|\Rightarrow w(E)\leq \delta w(Q)$ holds, whenever $Q$ is a d-dimensional cube and $E$ is its arbitrary measurable subset of $Q$. Note that $w$ satisfies the above condition if and only if it belongs $A_{p}$ class for some
	$p\in (1,\infty)$. Coifman and Fefferman \cite{CF} proposed another approach based on verifying the following inequality
	$$
	\frac{w(E)}{w(Q)}\leq C\left(\frac{|E|}{|Q|}\right)^{\varepsilon},
	$$
	where $Q,E$ are as before, while $C,\delta>0$ are constants depending only $w.$ Note that the two conditions lead to same class of weights. For More detailed information we refer the reader to \cite{DMO}.
	
	To prove the main result we need analogous result for Vilenkin group. It should be noted that we give the proof which we had not found in literature. 
	
	\begin{prop}\label{A_Infty Equivalent}
		Let $w\in A_{\infty}(I_{0})$, where $I_{0}\in \mathcal{F}$.  There exist positive constants $C,\varepsilon>0$ such that for any  generalized interval $I\subset I_{0}$ and measurable subset $E\subset I,$
		\begin{equation}\label{2.1}
			\frac{w(E)}{w(I)}\leq C\left(\frac{\mu(E)}{\mu(I)}\right)^{\delta}.
		\end{equation}
	\end{prop}
	
	For proving the result we need modified form of the Calderón-Zygmund decomposition lemma (see \cite{You1}, Lemma 2).
	
	\begin{lem}
		Given an interval $I\in \mathcal{F}$ and a function $f\in L^{1}(G)$, then for $t\geq|f|_{I}$, there exists  a collection
		$I_{j}$ of disjoint generalized intervals $I_{j}\subset I$  such that
		$$
		t<\frac{1}{\mu(I_{j})}\int_{I_{j}}|f|d\mu\leq3t,\,\,\,\forall I_{j},
		$$
		and for almost every $x\in I\backslash\cup_{j}I_{j},\,\,|f(x)|\leq t.$
	\end{lem}
	
	\emph{Proof of Proposition} \ref{A_Infty Equivalent}.
	Fix a generalized interval $I\subset I_{0}$ and for integer $k\geq 0$ define the sequence $t_{k}=10^{k}w_{I}=10^{k}t_{0}.$ Using Lemma 2.5  For each $k$  we may find Calderón-Zygmund generalized intervals $I_{j}^{k}$ of $w$ in following manner.  First construct Calderón-Zygmund generalized intervals $I_{j}^{0}$  relative to $I$ at height $t_{0}$ (Calderón-Zygmund generalized intervals of rang $0$). Denote $\Omega_{0}=\cup I^{0}_{j}$.  For any  fixed $I^{0}_{j}$ interval find  Calderón-Zygmund generalized intervals (of rang 1) of $w$ and height $t_{1}$. Denote by $I^{1}_{j}$ the intervals of rang $1$  and $\Omega_{1}=\cup I^{1}_{j}$. Note that $\Omega_{1}\subset\Omega_{0}\subset I.$ In this manner we may construct collection $I^{k}_{j}$ Calderón-Zygmund generalized intervals and the set $\Omega_{k}$ with properties:
	
	a) $\Omega_{k+1}\subset\Omega_{k},\,k=0,1,2,...$,
	
	b) $t_{k}< w_{I^{k}_{j}}\leq 3t_{k}, \,k=0,1,2,...$,
	
	c) $w(x)\leq t_{k}$, $x\in I\backslash \Omega_{k}.$
	
	Note that from the construction for any $i$ there exists $j$ such that $I^{k+1}_{i}\subset I^{k}_{j}.$
	
	Then
	$$
	\mu(\Omega_{k+1}\cap I_{j}^{k})=\sum_{I^{k+1}_{i}\subset I_{j}^{k}}\mu(I^{k+1}_{i})
	$$
	
	$$
	<t_{k+1}^{-1}\sum_{I^{k+1}_{i}\subset I_{j}^{k}}w(I^{k+1}_{i})\leq t_{k+1}^{-1}w(I^{k}_{j})\leq\frac{3t_{k}}{t_{k+1}}\mu(I^{k}_{j})=\frac{3}{10}\mu(I^{k}_{j}).
	$$
	
	Hence, by $A_{\infty}(I_{0})$ condition with $\varepsilon=3/10,$ there exists $\delta>0$ such that $w(\Omega_{k+1}\cap I_{j}^{k})\leq\delta w(I_{j}^{k}),$ 
	and if we sum over all $j$, we obtain $w(\Omega_{k+1})\leq\delta w(\Omega_{k})$ and consequently we have that $w(\Omega_{k})\leq \delta^{k+1}w(I).$
	
	For almost every $x\in I\setminus\Omega_{k},$  $w(x)\leq t_{k}$.  For fixed $\varepsilon$
	$$
	\frac{1}{\mu(I)}\int_{I}w(x)^{1+\varepsilon}d\mu=\frac{1}{\mu(I)}
	\int_{I\setminus\Omega_{0}}w(x)^{1+\varepsilon}d\mu+\frac{1}{\mu(I)}\sum_{k=0}^{\infty}\int_{\Omega_{k}\setminus\Omega_{k+1}}w(x)^{1+\varepsilon}d\mu
	$$
	
	$$
	\leq \frac{t_{0}^{\varepsilon}}{\mu(I)}
	\int_{I\setminus\Omega_{0}}w(x)d\mu+\frac{1}{\mu(I)}\sum_{k=0}^{\infty}t_{k+1}^{\varepsilon}w(\Omega_{k})
	$$
	
	$$
	\leq\frac{t_{0}^{\varepsilon}}{\mu(I)}
	\int_{I\setminus\Omega_{0}}w(x)d\mu+\frac{1}{\mu(I)}\sum_{k=0}^{\infty}10^{(k+1)\varepsilon}t_{0}^{\varepsilon}\delta^{k+1}w(I)
	$$
	
	Fix $\varepsilon>0$ so that $10^{\varepsilon}\delta<1$, we obtain that last term is bounded by
	$$
	t_{0}^{\varepsilon}\frac{1}{\mu(I)}\int_{I}w(x)d\mu+C\mu(I)^{-1}t_{0}^{\varepsilon}w(I)\leq C\left(\frac{1}{\mu(I)}\int_{I}w(x)d\mu\right)^{1+\varepsilon}.
	$$
	Hence, given $\varepsilon>0$ the weight satisfies Reverse H\"{o}lder inequality .
	
	Finally if we use H\"{o}lder's inequality for $w(E)=\int_E w(x)d\mu$ and Reverse H\"{o}lder's inequality for $1+\varepsilon$ we get \eqref{2.1}.\,\,
	\qed
	
	For $0<r<\infty$  define $M_{r}f(x)=M(|f|^{r})(x)^{1/r}.$  For brevity, hereafter we will write $f_{I}$ instead of $\int_{I}fd\mu/\mu(I).$
	
	As a consequence of the reverse H\"{o}lder inequality we get that if $w\in A_{p}(G)$ for some $p$, then there exists $s>1$ such that $M_{s}w(x)\leq C Mw(x).$  We need a sharper version of this inequality.
	
	\begin{prop}\label{Sharper Version}
		Given $w\in A_{1}(G)$, if $s_{0}=1+\frac{1}{8[w]_{A_{1}}}$, then for $1<s\leq s_{0}$ and for almost every $x,$
		\begin{equation}\label{Sharp ineq}
			M_{s}w(x)\leq 4 Mw(x)\leq4[w]_{A_{1}}w(x).
		\end{equation}
	\end{prop}
	This type of estimates is well known in Euclidian setting. For the sake of completeness  we will give a proof  for the Vilenkin group.

	We need an inequality that is the reverse of the weak $(1,1)$  inequality for maximal operator $M.$
	\begin{lem}\label{Reverse waek}
		Given a function $f\in L^{1}(G),$ for every interval $I\in\mathcal{F}$ and $t\geq|f|_{I},$
		$$
		\mu(\{x\in I\,:\,\,Mf(x)>t\})\geq\frac{1}{3t}\int_{\{x\in I:\,|f(x)|>t\}}|f(x)|d\mu.
		$$
	\end{lem}
	
	\emph{Proof.} $t\geq |f|_{I}$; if $t\geq\|f\|_{L^{\infty}},$ then this result is true. Otherwise, by Lemma 2.5, let $I_{i}$ be the Calderón-Zygmund intervals
	of $f$ relative to $I$ and $t.$ For every $x\in I_{i}$
	$$
	Mf(x)\geq\frac{1}{\mu(I_{i})}\int_{I_{i}}|f|d\mu>t.
	$$
	Since $|f(x)|\leq t$ for almost every $x\in I\backslash\cup_{i}I_{i},$ we have
	$$
	\mu(\{x\in I\,:\,\,Mf(x)>t\})\geq\sum_{j}\mu(I_{j})
	$$
	$$
	\geq\frac{1}{3t}\sum_{j}\int_{I_{j}}|f|d\mu\geq \frac{1}{3t}\int_{\{x\in I:\,|f(x)|>t\}}|f(x)|d\mu.
	$$
	\qed

	
	\emph{Proof of Proposition \ref{Sharper Version}}.
	Let $\varepsilon=(8[w]_{A_{1}})^{-1},\,\,s_{0}=1+\varepsilon,$ and fix an interval $I$ and $x_{0}\in I.$
	To prove the first inequality of \eqref{Sharp ineq} it is  sufficient to show that
	$$
	\frac{1}{\mu(I)}\int_{I}w(x)^{s_{0}}d\mu\leq 4Mw(x_{0})^{s_{0}}.
	$$
	
	We have that
	
	$$
	\frac{1}{\mu(I)}\int_{I}w(x)^{s_{0}}d\mu=\frac{1}{\mu(I)}\int_{I}w(x)^{\varepsilon}w(x)d\mu
	$$
	$$
	=\varepsilon(\mu(I))^{-1}\int_{0}^{\infty}t^{\varepsilon-1}w(\{x\in I:\,w(x)>t\})dt
	$$
	$$
	=\varepsilon(\mu(I))^{-1}\int_{0}^{Mw(x_{0})}t^{\varepsilon-1}w(\{x\in I:\,w(x)>t\})dt
	$$
	$$
	+\varepsilon(\mu(I))^{-1}\int_{Mw(x_{0})}^{\infty} t^{\varepsilon-1}w(\{x\in I:\,w(x)>t\})dt.
	$$
	
	For the first term we have
	$$
	\varepsilon(\mu(I))^{-1}\int_{0}^{Mw(x_{0})}t^{\varepsilon-1}w(\{x\in I:\,w(x)>t\})dt
	$$
	$$
	\leq\varepsilon(\mu(I))^{-1} w(I) \int_{0}^{Mw(x_{0})}t^{\varepsilon-1}dt=\frac{1}{\mu(I)}\int_{I}w(y)d\mu\cdot Mw(x_{0})^{\varepsilon}\leq Mw(x_{0})^{1+\varepsilon}.
	$$
	Using Lemma \ref{Reverse waek} we obtain
	$$
	\varepsilon(\mu(I))^{-1}\int_{Mw(x_{0})}^{\infty} t^{\varepsilon-1}w(\{x\in I:\,w(x)>t\})dt
	$$
	$$
	=\varepsilon(\mu(I))^{-1}\int_{Mw(x_{0})}^{\infty}t^{\varepsilon-1}\int_{\{x\in I:\,w(x)>t\}}d\mu dt
	$$
	$$
	\leq3\varepsilon(\mu(I))^{-1}\int_{0}^{\infty}t^{\varepsilon}\mu(\{x\in I;\,Mw(x)>t\})dt
	$$
	$$
	=\frac{3\varepsilon}{1+\varepsilon}\frac{1}{\mu(I)}\int_{I}Mw(x)^{1+\varepsilon}d\mu
	$$
	$$
	\leq\frac{3\varepsilon[w]_{A_{1}}^{1+\varepsilon}}{1+\varepsilon}\frac{1}{\mu(I)}\int_{I}w(x)^{1+\varepsilon}d\mu.
	$$

	From above estimates we get
	
	$$
	\frac{1}{\mu(I)}\int_{I}w(x)^{1+\varepsilon}d\mu\leq Mw(x_{0})^{1+\varepsilon}+\frac{3\varepsilon[w]_{A_{1}}^{1+\varepsilon}}{1+\varepsilon}\frac{1}{\mu(I)}\int_{I}w(x)^{1+\varepsilon}d\mu
	.$$
	
	Since for all $x\geq1, x^{1/8x}\leq2,$  we have
	$$
	\frac{3\varepsilon[w]_{A_{1}}^{1+\varepsilon}}{1+\varepsilon}\leq \frac{3}{8}[w]_{A_{1}}^{-1}[w]_{A_{1}}^{1+(8[w]_{A_{1}})^{-1}}\leq \frac{3}{4}
	$$
	
	and consequently the first inequality in \eqref{Sharp ineq} is valid. The second inequality in \eqref{Sharp ineq} is clear.
	\qed

	\section{Proof of Theorem \ref{Muck}}
	
	Given a generalized interval $I\in \mathcal{F}$ define the averaging operator $A_{I}$ by
	$$
	A_{I}f(x)=\frac{1}{\mu(I)}\int_{I}fd\mu\,\chi_{I}(x).
	$$
	
	\begin{prop}
		Given a exponent $p(\cdot),$ $1<p_{-}\leq p_{+}<\infty$, there exists a constant $C>0$ such that for any interval $I\in \mathcal{F}$
		$$
		\|A_{I}f\|_{p(\cdot)}\leq C\|f\|_{p(\cdot)}
		$$
		if and only if $p(\cdot)\in \mathcal{A}(G).$
	\end{prop}
	
	The proof of Proposition 3.1 is essentially the same as for averaging operator defined by cubes for Euclidean setting (see for example \cite{CUF}, Proposition 4.47).
	
	Lemma 3.2 shows that the condition $p(\cdot)\in \mathcal{A}(G)$ is actually sufficient for modular inequality.  Analogous estimate for the case $L^{p(\cdot)}(\mathbb{R}^{n})$ was obtained by Kopaliani \cite{K1}. The proof in \cite{K1} is based on some concepts from convex analysis.  Lerner  in \cite{L1} gave a different and simple proof. In this paper our approach is based on the adaptation of Lerner's proof \cite{L1}.
	
	\begin{lem}\label{main lemma}
		Given exponent $p(\cdot)$   such that $ 1<p_{-}\leq p_{+}<\infty,$ suppose  $p(\cdot)\in \mathcal{A}(G)$. Let $f\in L^{p(\cdot)}(G).$  If there exists an interval $I\in\mathcal{F}$ and constants $c_{1},c_{2}>0$ such that  $|f|_{I}\geq c_{1}$ and $\|f\|_{p(\cdot)}\leq c_{2},$ where $c_{1},c_{2}>0$, then there exists a constant $c$ depending only on $p(\cdot), c_{1},c_{2}$ such that
		\[
		\int_{I}(|f|_{I})^{p(x)}d\mu\leq c\int_{I}|f(x)|^{p(x)}d\mu.
		\]
	\end{lem}

	\textit{Proof.} Using the condition $p_{+}<\infty$  we may consider only the case $c_{1}=c_{2}=1.$ Since $p'_{+}<\infty$, there exists $\alpha>0$ such that
	\begin{equation}\label{3.2}
		\int_{I}\alpha^{p'(y)-1}d\mu=\int_{Q}|f(x)|d\mu.
	\end{equation}
	
	Since $|f|_{I}\geq1,$ we have $\alpha\geq1.$ By generalized H\"{o}lder inequality
	$$
	\int_{I}f(x)d\mu\leq 2\|f\|_{p(\cdot)}\|\chi_{I}\|_{p'(\cdot)}
	$$
	we get $\int_{I}\alpha^{p'(y)-1}d\mu\leq 2\|\chi_{I}\|_{p'(\cdot)}$ and consequently,
	\begin{equation}\label{3.4}
		\alpha\leq c/\|\chi_{I}\|_{p'(\cdot)}.
	\end{equation}
	
	Given this value $\alpha$, we have that
	\begin{equation}\label{3.5}
		\int_{I}(|f|_{I})^{p(x)}d\mu=\int_{I}\left(\frac{1}{\mu(I)}\int_{I}\alpha^{p'(y)-1}d\mu\right)^{p(x)}d\mu
	\end{equation}
	$$
	=\left(\frac{1}{\mu(I)}\int_{I}\left(\frac{1}{\mu(I)}\int_{I}\alpha^{p'(y)-p'(x)}d\mu\right)^{p(x)-1}d\mu\right)\int_{I}\alpha^{p'(y)}d\mu.
	$$
	
	For each $x\in I$ partition $I$ into
	$E_{1}(x)=\{y\in I:\, p'(y)>p'(x)\}$ and $E_{2}(x)=I\backslash E_{1}(x)$. Using \eqref{3.4} and the estimate  $\alpha\geq1$, we obtain
	$$
	\int_{I}\alpha^{p'(y)-p'(x)}d\mu=\int_{E_{1}(x)}\alpha^{p'(y)-p'(x)}d\mu+\int_{E_{2}(x)}\alpha^{p'(y)-p'(x)}d\mu
	$$
	$$
	\leq c(\|\chi_{I}\|_{p'(\cdot)})^{p'(x)}+\mu(I).
	$$
	In view of  $p(\cdot)\in A(G)$, we have
	\begin{equation}\label{3.6}
		\frac{1}{\mu(I)}\int_{I}\left(\frac{1}{\mu(I)}\int_{I}\alpha^{p'(y)-p'(x)}d\mu\right)^{p(x)-1}d\mu
	\end{equation}
	$$
	\leq c\frac{1}{\mu(I)}\int_{I}\left(\frac{1}{\mu(I)}(\|\chi_{I}\|_{p'(\cdot)})^{p'(x)}+1\right)^{p(x)-1}d\mu
	$$
	$$
	\leq c+c\frac{1}{\mu(I)}\int_{I}\left(\frac{1}{\mu(I)}(\|\chi_{I}\|_{p'(\cdot)})^{p'(x)}\right)^{p(x)-1}d\mu
	$$
	$$
	\leq c+c\int_{I}\left(\frac{\|\chi_{I}\|_{p'(\cdot)}}{\mu(I)}\right)^{p(x)}d\mu
	$$
	$$
	\leq c+c\int_{I}\left(\frac{1}{\|\chi_{I}\|_{p(\cdot)}}\right)^{p(x)}d\mu\leq c.
	$$
	Further,
	\begin{equation}\label{3.7}
		\int_{I}\alpha^{p'(y)}d\mu=2\alpha\int_{I}|f(x)|d\mu-\int_{I}\alpha^{p'(y)}d\mu
	\end{equation}
	$$
	\leq 2\alpha\int_{\{y\in I:\,2\alpha|f(y)|>\alpha^{p'(y)}\}}|f(y)|d\mu
	$$
	$$
	\leq c\int_{I}|f(y)|^{p(y)}d\mu.
	$$
	
	From \eqref{3.5}, \eqref{3.6} and \eqref{3.7} we obtain desired estimate.    \qed
	
	\begin{cor}\label{Ainfinity}
		Let $1<p_{-}\leq p_{+}<\infty$ and $p(\cdot)\in\mathcal{A}(G).$ Suppose that $\xi_{1}\leq t\leq \xi_{2}/\|\chi_{I}\|_{p(\cdot)}$, where $\xi_{1},\xi_{2}>0$ and $I\in \mathcal{F}.$ Then $t^{p(x)}\in A_{\infty}(I)$ with $A_{\infty}$ constant depending only on $p(\cdot),\,\xi_{1},\xi_{2}.$
	\end{cor}

	\textit{Proof.} Let $I'\subset I$, where $I', I \in\mathcal{F}$ and $E\subset I'$ be any measurable subset with  $\mu(E)>\mu(I')/2.$ Define $f=t\chi_{E}.$ Then
	$$
	|f|_{I'}=\frac{1}{\mu(I')}\int_{I'}t\chi_{E}(x)d\mu=t\frac{\mu(E)}{\mu(I')}\geq\frac{\xi_{1}}{2},
	$$
	$$
	\|f\|_{p(\cdot)}=t\|\chi_{E}\|_{p(\cdot)}\leq\xi_{2}\frac{\|\chi_{E}\|_{p(\cdot)}}{\|\chi_{I}\|_{p(\cdot)}}\leq\xi_{2}.
	$$
	
	Therefore, $f$ satisfies the hypotheses of  Lemma 3.2 with $c_{1}=\xi_{1}/2,\,c_{2}=\xi_{2}$ and there exists a constant $c$ depending only on $p(\cdot), \xi_{1},\xi_{2}$ such that
	$$
	\frac{1}{2^{p_{+}}}\int_{I_{0}}t^{p(\cdot)}d\mu\leq c\int_{E}t^{p(\cdot)}d\mu,
	$$
	which proves that  $t^{p(x)}\in A_{\infty}(I).$ \qed
	
	\textit{Proof of Theorem \ref{Muck}}. The part $(ii)\Rightarrow (i)$ of Theorem 1.6 follows immediately from Proposition 3.1 and from the fact that $|f|_{I}\chi_{I}(x)\leq Mf(x)$ for any interval $I\in\mathcal{F}.$
	
	Implication $(i)\Rightarrow(ii).$
	Suppose $f\in L^{p(\cdot)}(G)$ and $\|f\|_{p(\cdot)}\leq1$. It is sufficient to proof that there exists a positive constant $C$ (independent of $f$) such that for any nonnegative function  $g\in L^{p'(\cdot)}(G)$, with  $\|g\|_{p'(\cdot)}\leq1$
	\begin{equation}\label{3.8}
		\int_{G}Mf(x)g(x)d\mu\leq C.
	\end{equation}
	
	For each positive integer $k$ set
	$$
	\Omega_{k}=\{x\in G\,\,:\,\,Mf(x)>3^{k}\}.
	$$
	
	Note that
	\begin{equation}\label{3.88}
		\int_{G\backslash\Omega_{1}}Mf(x)g(x)d\mu\leq C.
	\end{equation}
	
	Define $D_{k}=\Omega_{k}\backslash\Omega_{k+1}.$ Let $F_{k}$ be an arbitrary compact subset of $D_{k}.$ We will prove that
	\begin{equation}\label{3.9}
		\int_{\cup F_{k}}Mf(x)g(x)d\mu\leq C.
	\end{equation}
	By simple limiting argument from \eqref{3.9} and from \eqref{3.88} we obtain \eqref{3.8}.
	
	Let $\mu(F_{k})>0.$  There exists a finite collection of generalized intervals $I_{\alpha},\alpha\in A_{k} $, $F_{k}\subset\cup_{\alpha\in A_{k}}I_{\alpha}$, such that $|f|_{I_{\alpha}}>3^{k},\,\,\alpha\in A_{k}$ and for all fixed $\alpha$, there exists $x_{\alpha}\in I_{\alpha}$ such that $ Mf(x_{\alpha})\leq 3^{k+1}.$  Note that if $I_{\alpha_{1}}$ and $I_{\alpha_{2}}$ belong to distinct $\mathcal{F}_{l}$'s and are not disjoint ($\mu(I_{\alpha_{1}}\cap I_{\alpha_{2}})>0$) then one is a subset of the other. Consequently  without loss  of generality we may assume that in collection $I_{\alpha},\alpha\in A_{k} $ if $\mu(I_{\alpha_{1}}\cap I_{\alpha_{2}})>0$ for some  $\alpha_{1}$ and $\alpha_{2}$, then $I_{\alpha_{1}}$ and $I_{\alpha_{2}}$ belong to the same $\mathcal{F}_{l}$'s (for some $l$). By Vitali covering lemma, we may select from collection $I_{\alpha},\alpha\in A_{k}$ the finite collection of pairwise disjoint intervals $\{I_{j}^{k}\}$ $j\in\{1,...,N_{k}\}$ such that $F_{k}\subset\cup_{j}3I_{j}^{k}.$
	
	Without loss of generality we may assume that  $\mu(F_{k})>0$ for all $k\geq1.$ Define the sets $E_{1}^{k}=3I_{1}^{k}\cap F_{k,}$ $E_{j}^{k}=(3I_{j}^{k}\backslash\cup_{s<j}3I_{s}^{k})\cap F_{k},\,j>1.$ Note that the sets $E_{j}^{k}$ are pairwise disjoint and $\cup_{j}E_{j}^{k}=F_{k}.$

	Define
	\[
	Tg(x)=\sum_{k=1}^{\infty}\sum_{j}\left(\frac{1}{\mu(I_{j}^{k})}\int_{E_{j}^{k}}gd\mu\right)\chi_{I_j^{k}}(x).
	\]
	
	Using the above definition, we get
	\[
	\int_{\cup_{k}F_{k}}(Mf)(x)g(x)d\mu\leq3^{k+1}\sum_{k=1}^{\infty}\sum_{j}\int_{E_{j}^{k}}gd\mu \leq 3\sum_{k=1}^{\infty}\sum_{j}f_{I_{j}^{k}}\int_{E_{j}^{k}}gd\mu
	\]
	\[
	=3\int_{G}fTg\leq
	6\|f\|_{p(\cdot)}\|Tg\|_{p'(\cdot)},
	\]
	and consequently for proving\eqref{3.9}, it is sufficient to show that $\|Tg\|_{p'(\cdot)}\leq C.$

	Note that $I_{j}^{k}\subset\Omega_{k}=\cup_{l=0}^{\infty}D_{k+l}$ and hence $Tg=\sum_{l=0}^{\infty}T_{l}g,$ where
	\[
	T_{l}g(x)=\sum_{k=1}^{\infty}\sum_{j}a_{j,k}(g)\chi_{I_{j}^{k}\cap D_{k+l}}(x),\,\,\,(l=0,1,...)
	\]
	where  $\alpha_{j,k}(g)=\frac{1}{\mu(I_{j}^{k})}\int_{E_{j}^{k}}gd\mu.$
	
	Let $\mathcal{I}_{1}=\{(j,k)\,\,:\,\,\alpha_{j,k}(g)>1\}$ and $\mathcal{I}_{2}=\{(j,k)\,\,:\,\,\alpha_{j,k}(g)\leq1\}.$
	
	By condition $p\in \mathcal{A}(G)$ and H\"{o}lder inequality implies that for any interval $I\in\mathcal{F},$ $\|\chi_{3I}\|_{p(\cdot)}\leq C\|\chi_{I}\|_{p(\cdot)}.$ We have
	
	\[
	\alpha_{j,k}(g)\leq\frac{2}{\mu(I_{j}^{k})}\|\chi_{E_{j}^{k}}\|_{p(\cdot)}\|g\chi_{E_{j}^{k}}\|_{p'(\cdot)}
	\leq\frac{2}{\mu(I_{j}^{k})}\|\chi_{3I_{j}^{k}}\|_{p(\cdot)}
	\]
	\[
	\leq\frac{C}{\|\chi_{3I_{j}^{k}}\|_{p'(\cdot)}}\leq\frac{C}{\|\chi_{I_{j}^{k}}\|_{p'(\cdot)}}.
	\]
	
	Let $(j,k)\in \mathcal{I}_{1}.$ Then by Corollary \ref{Ainfinity} $\alpha_{j,k}(g)^{p'(x)}\in A_{\infty}(I_{j}^{k})$
	and by Lemma \ref{main lemma}, (see, also \eqref{2.1})
	\[
	\int_{I_{j}^{k}\cap D_{k+l}}\alpha_{j,k}(g)^{p'(x)}d\mu\leq C\left(\frac{\mu(I_{j}^{k}\cap D_{k+l})}{\mu(I_{j}^{k})}\right)^{\varepsilon}\int_{I_{j}^{k}}\alpha_{j,k}(g)^{p'(x)}d\mu
	\]
	\begin{equation}\label{3.10}
		\leq C\left(\frac{\mu(I_{j}^{k}\cap D_{k+l})}{\mu(I_{j}^{k})}\right)^{\varepsilon}\int_{E_{j}^{k}}g(x)^{p'(x)}d\mu.
	\end{equation}
	
	If $(j,k)\in \mathcal{I}_{2}$, then we have
	$$
	\int_{I_{j}^{k}\cap D_{k+l}}\alpha_{j,k}(g)^{p'(x)}d\mu\leq\int_{I_{j}^{k}\cap D_{k+l}}\alpha_{j,k}(g)d\mu
	$$
	\begin{equation}\label{3.11}
		=\frac{\mu(I_{j}^{k}\cap D_{k+l})}{\mu(I_{j}^{k})}\int_{E_{j}^{k}}g(x)d\mu.
	\end{equation}
	
	We need estimate  $\mu(I_{j}^{k}\cap D_{k+l})$ for $l\geq2.$ Let $x\in I_{j}^{k}$ and $I\in\mathcal{F}$ be an arbitrary interval such that $x\in I.$ Observe that either $I\subset 3I_{j}^{k}$ or $I_{j}^{k}\subset3I.$ If the second inclusion holds, then $3I\cap D_{k}\neq\emptyset$ and hence
	$$
	|f|_{I}\leq 3|f|_{3I}\leq3\cdot3^{k+1}\leq3^{k+l}\,\,(l\geq2).
	$$
	Therefore, if $|f|_{I}>3^{k+l}$, then $I\subset3I_{j}^{k}.$ From this and from weak type property of $M,$ we get
	$$
	\mu(I_{j}^{k}\cap D_{k+l})\leq\mu\{x\in I_{j}^{k}:\,\,M(f\chi_{3I_{j}^{k}})(x)>3^{k+l}\}
	$$
	\begin{equation}\label{3.12}
		\leq\frac{C}{3^{k+l}}\int_{3I_{j}^{k}}|f|d\mu\leq C\frac{\mu(I_{j}^{k})}{3^{k+l}}|f|_{3I_{j}^{k}}\leq C\frac{3^{k+1}}{3^{k+l}}\mu(I_{j}^{k})\leq\frac{C}{3^{l}}\mu(I_{j}^{k}).
	\end{equation}

	By estimates \eqref{3.10},\eqref{3.11},\eqref{3.12}, when $l\geq2$ we obtain
	\[
	\int_{G}(T_{l}g(x))^{p'(x)}d\mu=\sum_{k=1}^{\infty}\sum_{j}\int_{I_{j}^{k}\cap D_{k+l}}\alpha_{j,k}(g)^{p'(x)}d\mu
	\]
	\[
	\leq C3^{-l\varepsilon}\sum_{(j,k)\in\mathcal{I}_{1}}\int_{E_{j}^{k}}g(x)^{p'(x)}d\mu+C3^{-l}\sum_{(j,k)\in\mathcal{I}_{2}}\int_{E_{j}^{k}}g(x)d\mu
	\]
	\[
	\leq C3^{-l\alpha}\left(\int_{G}g(x)^{p'(x)}d\mu+\int_{G}g(x)d\mu\right).
	\]
	Where $\alpha=\min\{ {1}, {\varepsilon} \}$.
	
	Using the fact that $\|g\|_{1}\leq 2\|\chi_{G)}\|_{p'(\cdot)}$, and $\int_{G}g(x)^{p'(x)}d\mu\leq1$ we obtain
	\[
	\|T_{l}g\|_{p'(\cdot)}\leq C3^{-l\alpha/p'_{+}}\,\,\,(l\geq2).
	\]
	
	For $l=0,1$  if we  use a trivial estimate $\mu(I_{j}^{k}\cap D_{k+l})\leq\mu(I_{j}^{k})$, analogously will be obtained the estimate $\|T_{l}g\|_{p'(\cdot)}\leq C$. Finally we obtain
	\[
	\|Tg\|_{p'(\cdot)}\leq\sum_{l=0}^{\infty}\|T_{l}g\|_{p'(\cdot)}\leq C.
	\]
	\qed

	\section{Proof of Theorem 1.8}

	The implication $(ii)\Rightarrow(i)$ is straightforward. Fix $r_{0},\,r_{0}<r<1,$ and let $s=1/r.$ by H\"{o}lder's inequality, we have that $Mf(x)\leq M(|f|^{s})(x)^{1/s}=M_{s}f(x).$  Note that $\||f|^{s}\|_{p(\cdot)}=\|f\|_{sp(\cdot)}^{s}$  and
	$$
	\|Mf\|_{p(\cdot)}\leq\|M(|f|^{s})^{1/s}\|_{p(\cdot)}=\|M(|f|^{s})\|_{rp(\cdot)}^{r} \leq C  \||f|^{s}\|_{rp(\cdot)}^{r}=C\|f\|_{p(\cdot)}.
	$$
	
	To prove that $(i)\Rightarrow(ii),$ we first construct a $A_{1}(G)$ weight using the Rubio de Francia iteration algorithm.
	Given $h\in L^{p(\cdot)}(G),$ define
	$$
	\mathcal{R}h(x)=\sum_{k=0}^{\infty}\frac{M^{k}h(x)}{2^{k}\|M\|_{L^{p(\cdot)}(G)}^{k}},
	$$
	where for $k\geq1,\,\,M^{k}=M\circ M\circ\cdot\cdot\cdot\circ M$ denotes $k$ iterations of the Maximal operator $M$ and $M^{0}f=|f|.$ The function $\mathcal{R}h(x)$
	has the following properties:
	
	(a) For all $x\in G,$ $|h(x)|\leq \mathcal{R}h(x);$
	
	(b) $\mathcal{R}$ is bounded on $L^{p(\cdot)}(G)$ and $\|\mathcal{R}h\|_{p(\cdot)}\leq 2\|h\|_{p(\cdot)};$
	
	(c) $\mathcal{R}h\in A_{1}(G)$ and $[\mathcal{R}h]_{A_{1}}\leq 2\|M\|_{L^{p(\cdot)}(G)}.$
	
	The proof of properties (a),(b),(c) are the same, as Euclidian setting (see \cite{CUF}, pp.157) and we omit it here. By property (c) and Proposition \ref{Sharper Version} there exists $s_{0}>1$ such that for all $s,\,1<s<s_{0},$
	$$
	M_{s}(\mathcal{R}h)(x)\leq M_{s_{0}}(\mathcal{R}h)(x)\leq 8 \|M\|_{L^{p(\cdot)}(G)}\mathcal{R}h(x).
	$$
	
	Let $r_{0}=1/s_{0}$. Fix $r$ such that $r_{0}<r<1.$ Let $s=1/r.$

	By properties (a) and (b) we have
	$$
	\|Mf\|_{rp(\cdot)}= \|(Mf)^{1/s}\|_{p(\cdot)}^{s}= \|M_{s}(|f|^{r})\|_{p(\cdot)}^{s}\leq \|M_{s}(\mathcal{R}(|f|^{r})\|_{p(\cdot)}^{s}
	$$
	$$
	\leq C\|M\|_{L^{p(\cdot)}(G)}^{s}\|\mathcal{R}(|f|^{r})\|_{p(\cdot)}^{s}\leq C\||f|^{r}\|_{p(\cdot)}^{s}=C\|f\|_{rp(\cdot)}.
	$$
	\qed
	
	\section{Proof of Theorem 1.10}
	Since Vilenkin  polynomials are dense in $L^{p(\cdot)}(G)$ ($1\leq p_{-}\leq p_{+}<\infty$) the proof of equivalence of $(ii)$ and $(iii)$ is straightforward.
	The implications $(i)\Rightarrow(ii)$ follows from Rubio de Francia extrapolation theorem (Theorem 1.9), if we use Young's  weighted estimates for partial sum $S_{n}f$ of the Vilenkin-Fourier series (Theorem 1.2), Theorem 1.6, Theorem 1.8 and corollary 1.7.
	
	\textit{Proof of} $(ii)\Rightarrow(i).$ Consider $I\in \mathcal{F}.$  There is $x\in G$ such that $I$ is a proper subset of $x+G_{k}$ and $I$
	is a union of cosets of $G_{k+1}.$ First consider the case $\mu(I)\leq \mu(G_{k})/2.$
	Take
	$\alpha_{k}=[\mu(G_{k})/2\mu(I)]$, where $[a]$ is the largest integer less than or equal to $a$.  We have $\alpha_{k}\geq1.$ Let $f\in L^{p(\cdot)}(G)$ be a nonnegative function with support in $I.$  We use the following estimate
	(see \cite{You2}, pp.286-287): for $x\in I,$
	
	$$
	\phi_{k}^{-(\alpha_{k}-1)/2}(x)S_{\alpha_{k}m_{k}}(f\phi_{k}^{(\alpha_{k}-1)/2})(x)\geq\frac{1}{2\pi\mu(I)}\int_{I}f(t)d\mu=\frac{1}{2\pi}A_{I}f(x).
	$$
	We have
	$$
	\|A_{I}f\|_{p(\cdot)}\leq C\|\phi_{k}^{-(\alpha_{k}-1)/2}S_{\alpha_{k}m_{k}}(f\phi_{k}^{-(\alpha_{k}-1)/2})\|_{p(\cdot)}\leq C\|f\|_{p(\cdot)}.
	$$
	
	From this estimate  we obtain in standard way \eqref{vMuk}  in case $\mu(I)\leq \mu(G_{k})/2$ (see Proposition 3.1).
	
	Consider the case $\mu(I)> \mu(G_{k})/2.$ Note that every coset of $G_{k}$ is in $\mathcal{F}_{k-1}$ and $\mu(G_{k})\leq \mu(G_{k-1})/2$ and  consequently \eqref{vMuk} holds for all cosets of $G_{k}$. We have
	$$
	\|\chi_{I}\|_{p(\cdot)}\|\chi_{I}\|_{p'(\cdot)}\leq \|\chi_{x+G_{k}}\|_{p(\cdot)}\|\chi_{x+G_{k}}\|_{p'(\cdot)}\leq C\mu(G_{k})\leq C\mu(I).
	$$
	\qed

\end{document}